\documentclass[11pt,a4paper,equation]{article}
\usepackage[greek,american]{babel}
\usepackage[iso-8859-7]{inputenc}
\usepackage{amssymb,amsfonts}
\usepackage[dvips]{graphicx}
\usepackage{amsmath}
\usepackage{epsfig}
\usepackage{latexsym}
\usepackage{makeidx}
\usepackage{pst-math,pst-xkey}

\numberwithin{equation}{section}
\makeindex
\begin{document}
\date{}
\author{Vassilis G. Papanicolaou \\
Department of Mathematics\\
National Technical University of Athens\\
Zografou Campus\\
157 80 Athens, GREECE\\
\underline{e-mail}: papanico@math.ntua.gr}
\title{An Arctangent Law}
\maketitle
\begin{abstract}
Let $M_r$ be the maximum value of an one-dimensional Brownian motion on the (time) interval $[0, r]$. We derive an explicit formula for the distribution of the time required (after $r$) for the Brownian motion to exceed $M_r$.
\end{abstract}
\textbf{2010 AMS Mathematics Classification.} 60J65.
\section{Introduction--The Theorem}
Let $B = \{B_t\}_{t \geq 0}$ be a standard Brownian motion in $\mathbb{R}^1$. The distribution of $B_0$ can be arbitrary. We set
\begin{equation*}
M_r := \max_{0 \leq t \leq r} B_t,
\end{equation*}
where $r > 0$ is a fixed time. The object of this brief note is the calculation of the distribution of the random time
\begin{equation}
S := \inf\{t \geq r : B_t = M_r\} - r = \inf\{t \geq r : B_t > M_r\} - r
\label{B0}
\end{equation}
(obviously, $S + r$ is a stopping time for $B$).\\\\
\textbf{Theorem.} \textit{The distribution function of the random variable $S$ defined in (\ref{B0}) is}
\begin{equation}
F_S(s) := P\{S \leq s\} = \frac{2}{\pi} \arctan\left(\sqrt{\frac{s}{r}} \, \right),
\qquad
s \geq 0.
\label{B00}
\end{equation}
It is remarkable that $F_S(s)$ is an elementary function.
\section{A Lemma}
Suppose that $X > 0$ is a random variable and $W = \{W_t\}_{t \geq 0}$ is an one-dimensional Brownian motion
with $W_0 = 0$. The passage time of $W$ to the level $X$ is
\begin{equation}
T_X := \inf\{t \geq 0 : W_t = X\}.
\label{A1}
\end{equation}
We remind the reader that if $X$ is not random, say $X = x > 0$, then the reflection principle (see, e.g., \cite{K-S}) yields
\begin{equation}
P\{T_x \leq t\} = 2 P\{W_t \geq x\} = \sqrt{\frac{2}{\pi}} \int_{x / \sqrt{t}}^\infty \; e^{-\xi^2 / 2} \, d\xi,
\qquad
t > 0.
\label{B1}
\end{equation}
\textbf{Lemma.} \textit{If $X$ and $W$ are independent, then the density of $T_X$ is}
\begin{equation}
f_{T_X}(t) = \frac{1}{\sqrt{2 \pi t^3}} \, E\left[ X e^{-X^2 / 2 t} \right],
\qquad
t > 0.
\label{B3}
\end{equation}
\textit{Proof}. For $t > 0$ we have
\begin{equation*}
P\{T_X \leq t\} = \int_0^\infty P\{T_X \leq t, X \in dx\} = \int_0^\infty P\{T_X \leq t \, |\, X = x\} \, dF_X(x)
\end{equation*}
\begin{equation*}
= \int_0^\infty P\{T_x \leq t \, |\, X = x\} \, dF_X(x) = \int_0^\infty P\{T_x \leq t\} \, dF_X(x),
\end{equation*}
where the last equality follows from the independence of $X$ and $W$. Thus, by invoking (\ref{B1})
\begin{equation}
P\{T_X \leq t\} = \sqrt{\frac{2}{\pi}} \int_0^\infty \left[\int_{x / \sqrt{t}}^\infty \; e^{-\xi^2 / 2} \, d\xi\right] dF_X(x),
\label{B2}
\end{equation}
from which (\ref{B3}) follows by differentiation with respect to $t$ (the passing of $d / dt$ inside the integral with respect to $x$
is justified by the fact that the quantity $x \exp(-x^2 / 2 t)$ is bounded in $x$, for any $t > 0$). \hfill $\blacksquare$\\\\
Before closing this short section it maybe worth noticing that, by applying integration by parts (or Tonelli) in the integral of (\ref{B2})
one obtains the formula
\begin{equation*}
P\{T_X \leq t\} = \sqrt{\frac{2}{\pi t}} \int_0^\infty e^{-x^2 / 2 t} F_X(x) \, dx,
\qquad
t > 0,
\end{equation*}
where $F_X(x)$ is the distribution function of $X$.
\section{Proof of the Theorem}
As before, let $r$ be a given positive number. We set
\begin{equation}
X = M_r - B_r \qquad\qquad \text{and} \qquad\qquad W_t = B_{t + r} - B_r, \qquad t \geq 0.
\label{B4}
\end{equation}
Clearly, $X > 0$ a.s., $W = \{W_t\}_{t \geq 0}$ is a Brownian motion with $W_0 = 0$ and,
also, $W$ and $X$ are independent. Furthermore, as it was observed by P. L\'{e}vy (see, e.g., \cite{K-S}, p. 97), $X$
and $|B_r - B_0|$ have the same law and, therefore, the density of $X$ is
\begin{equation}
f_X(x) = \sqrt{\frac{2}{\pi r}} \, e^{-x^2 / 2 r},
\qquad
x > 0.
\label{B5}
\end{equation}
We are now ready to complete the proof of our result.\\\\
\textit{Proof of the theorem}. Observe that, in view of (\ref{A1}) the random time $S$ of (\ref{B0}) can be expressed as
\begin{equation*}
S = T_X,
\end{equation*}
where $X$ is given by (\ref{B4}). Thus, we can use (\ref{B3}) to obtain the density of $S$
\begin{equation*}
f_S(s) = \frac{1}{\sqrt{2 \pi s^3}} \int_0^\infty x e^{-x^2 / 2 s} f_X(x) \, dx,
\qquad
s > 0,
\end{equation*}
where $f_X(x)$ is given by (\ref{B5}). It follows that
\begin{equation*}
f_S(s) = \frac{1}{\pi \sqrt{r s^3}} \int_0^\infty x e^{-x^2 / 2 s} e^{-x^2 / 2 r} \, dx
= \frac{\sqrt{r}}{\pi} \cdot \frac{1}{(s + r)\sqrt{s}} \, ,
\qquad
s > 0,
\end{equation*}
from which (\ref{B00}) follows immediately by integration with respect to $s$. \hfill $\blacksquare$\\\\
\textbf{Remark.} A remarkable consequence of the theorem is the following: For given $r_1$, $r_2$, with $0 \leq r_1 < r_2$, we set
\begin{equation*}
M_{[r_1, r_2]} := \max_{r_1 \leq t \leq r_2} B_t
\end{equation*}
and
\begin{equation*}
S_{[r_1, r_2]} := \inf\left\{t \geq r_2 : B_t = M_{[r_1, r_2]} \right\} - r_2.
\end{equation*}
Then,
\begin{equation*}
P\left\{S_{[r_1, r_2]} \leq s \right\} = \frac{2}{\pi} \arctan\left(\sqrt{\frac{s}{r_2 - r_1}} \, \right),
\qquad
s \geq 0.
\end{equation*}



\begin{thebibliography}{4}

\bibitem{K-S} I. Karatzas and S.E. Shreve, \textit{Brownian Motion and Stochastic Calculus}, Second
Edition, Springer, New York, 1991.

\end{thebibliography}
\end{document}